\documentclass[11pt,twoside]{article}

\newcommand{\sub}{\supseteq}

\def\ifif {if and only if\ \ }
\def\sub {\subseteq}

\usepackage{amsmath}
\usepackage{amsfonts}
\usepackage{amssymb}
\usepackage{makeidx}
\usepackage{amsthm}
\newcommand{\mathscr}{\mathbf}

\begin{document}
\title{\textsc{ some new
classes of topological spaces and annihilator ideals}}
%\dedication{\langle Professor Azarpanah \rangle}
%\author{ F. Azarpanah \\Department of Mathematics, Chamran University, Ahvaz , Iran}
%\thanks{2000 {\it Mathematics Subject Classification}: 54C40.}
\author{A. Taherifar \\ Departement of Mathematics, Yasouj University, Yasouj,
Iran}
\date{ataherifar@mail.yu.ac.ir\\
\[{\text{ Dedicated to \textsc{Professor} William Wistar
Comfort}}\]}

\maketitle
%\maketitle
%--------------------------------------------------------
\theoremstyle{definition}\newtheorem{thm}{Theorem}[section]
\theoremstyle{definition}\newtheorem{cor}[thm]{Corollary}
\theoremstyle{definition}\newtheorem{lem}[thm]{Lemma}
\theoremstyle{definition}\newtheorem{prop}[thm]{Proposition}
\theoremstyle{definition}\newtheorem{defn}[thm]{Definition}
\theoremstyle{definition}\newtheorem{Rem}[thm]{Remark}
\theoremstyle{remark}\newtheorem{rem}[thm]{Remark}
\theoremstyle{definition}\newtheorem{exam}[thm]{Example}
\theoremstyle{definition}\newtheorem{ques}[thm]{Questions}
%---------------------------------------------------------------------

{\noindent\bf Abstract.} By a characterization of semiprime
$SA$-rings by Birkenmeier, Ghirati and Taherifar in \cite[Theorem
4.4]{B}, and by the topological characterization of $C(X)$ as a
Baer-ring by Stone and Nakano in \cite[Theorem 3.25]{KM}, it is
easy to see that $C(X)$ is an $SA$-ring (resp., $IN$-ring) \ifif
$X$ is an extremally disconnected space. This result motivates the
following questions: Question $(1)$: What is $X$ if for any two
ideals $I$ and $J$ of $C(X)$ which are generated by two subsets of
idempotents, $Ann(I)+Ann(J)=Ann(I\cap J)?$ Question $(2)$: When
does for any ideal $I$ of $C(X)$ exists a subset $S$ of
idempotents such that $Ann(I)=Ann(S)$? Along the line of answering
these questions we introduce two classes of topological spaces. We
call $X$ an $\textit{EF}$ (resp., $\textit{EZ}$)-$\textit{space}$
if disjoint unions of clopen sets are completely separated (resp.,
every regular closed subset is the closure of a union of clopen
subsets). Topological properties of $\textit{EF}$ (resp.,
$\textit{EZ}$)-$\textit{spaces}$ are investigated. As a
consequence, a completely regular Hausdorff space $X$ is an
$F_{\alpha}$-space in the sense of Comfort and Negrepontis for
each infinite cardinal $\alpha$ \ifif $X$ is an $EF$ and
$EZ$-space. Among other things, for a reduced ring $R$ (resp.,
$J(R)=0$) we show that $Spec(R)$ (resp., $Max(R)$) is an
$EZ$-space \ifif for every ideal $I$ of $R$ there exists a subset
$S$ of idempotents of $R$ such that $Ann(I)=Ann(S)$.
%\end{abstract}
%--------------------------------------------------------------

\vspace{0.3 cm}\noindent{\it AMS Subject Classification Primary:
$54G05$, $54C40$. Secondary: 13A15.}

 \vspace{0.3 cm}\noindent{\it Keywords:}
$F_{\alpha}$-space, Extremally disconnected space,
Zero-dimensional space, $EF$-space, $EZ$-space, Reduced ring.

\section{Preliminaries}
A space $X$ is $\textit{extremally disconnected}$ (resp.,
$\textit{basically disconnected}$) if the closure of every open
subset is clopen in $X$ (resp., if the closure of any cozeroset is
open). It is well known that $X$ is an extremally disconnected
space \ifif any two disjoint open subsets of $X$ are completely
separated \ifif every open subset of $X$ is $C^{*}$-embedded. In
an extremally disconnected space all dense subsets are
$C^{*}$-embedded. The reader is referred to $\cite[1. H]{G}$,
\cite{H} and \cite{ST}. A topological space is said to be
$\textit{zero-dimensional}$  if it is a non-empty $T_{1}$-space
with a base consisting of clopen sets. Zero-dimensional spaces
were defined by Sierpi´nski in \cite{S}. All zero-dimensional
spaces are completely regular. A zero-dimensional space need not
be a normal space. The space $\beta T=W^{*}\times N^{*}$ where $T$
is the Tychonoff plank is an example of non-normal
zero-dimensional space (see \cite[Example 16.18] {G}). A space $X$
is $\textit{totally disconnected}$ \ifif the components in $X$ are
the singletons. Equivalently, $X$ is totally disconnected \ifif
the only non-empty connected subsets of $X$
 are the one-point sets. The following implications
characterize the relationship among the notions defined above: $X$
is extremally disconnected and $T_{3}$  $\Rightarrow$ $X$ is
zero-dimensional $\Rightarrow$ $X$ is totally disconnected. None
of the implications can be reversed and counterexamples exist even
in the class of metric spaces. For terminology and notations, the
reader is referred to \cite{E} and \cite{G}.
\\\indent All rings are assumed to be commutative with identity. By a reduced ring, we mean a ring without nonzero nilpotent elements.
For each subset $S$ of a ring $R$, $Ann(S)=\{r\in R: rs=0, \forall
s\in S\}$. A ring $R$ is called a $\textit{SA}$ (resp.,
$\textit{IN}$)-$\textit{ring}$ if for any two ideals $I, J$ of
$R$, $Ann(I)+Ann(J)=Ann(K)$, for some ideal $K$ of $R$ (resp.,
$Ann(I)+Ann(J)=Ann(I\cap J)$) (see \cite{B}). We denote the
Jacobson radical of $R$ by $J(R)$. For terminology and notations,
the reader is referred to \cite{HU}.\\\indent Throughout this
paper, we denote by $C(X)$, the ring of all real-valued continuous
functions on a completely regular Hausdorff space $X$, and
$C^{*}(X)$ is its subring of bounded functions. A completely
regular Hausdorff space $X$ is an $F$-$space$ if its cozerosets
are $C^{*}$-embedded. Equivalently, $X$ is an $F $-space if
finitely generated ideals of $C(X)$ are principal. For terminology
and notations, the reader is referred to \cite{G}.

For the  proof of the following lemma see \cite[ Corollary
3.6.5]{E}.
\begin{lem}\label{1.1}
If $A$ is a clopen subset of a topological space $X$, then
cl$_{\beta X}A$ is a clopen subset of $\beta X$.
\end{lem}
For the proof of the following theorems see \cite[ 1.17]{G} and
\cite[1.15]{G}.
\begin{thm}\label{1.2}
 A subset $S$ of $X$ is $C^{*}$-embedded in $X$ \ifif any two completely separated sets in $S$ are completely separated in $X$.
\end{thm}
\begin{thm}\label{1.3}
Two sets are completely separated \ifif they are contained in
disjoint zero-sets. Moreover, completely separated sets have
disjoint zero-set-neighborhoods.
\end{thm}
%\section{Preliminary Notes / Materials and Methods} Preliminary notes, materials and methods used in the paper.

\section{$EF$-space}
We call a topological space $X$ an $\textit{EF}$-$\textit{space}$
if for any two collections $\mathcal{U}$ and $\mathcal{V}$  of
clopen subsets of $X$ with $\bigcup \mathcal{U}\cap \bigcup
\mathcal{V}=\emptyset$, we have $\bigcup \mathcal{U}$ and $\bigcup
\mathcal{V}$ are completely separated. The class of $EF$-space
contains the class of spaces which are sums of connected spaces,
all spaces for which the closure of any union of clopen subsets is
open (hence all extremally disconnected spaces) and  all spaces
which any union of clopen subsets is $C^{*}$-embedded. If we take
$X$ as the sum of $\Bbb R$ (with usual topology) and $\Bbb N$,
then $X$ is an $EF$-space which is neither connected nor
extremally disconnected. A zero-dimensional space need not be an
$EF$-space. In fact, it is straightforward to check that if $X$ is
a zero-dimensional space, then $X$ is an $EF$-space if and only if
it is extremally disconnected. In this section, we prove that for
any two ideals $I$ and $J$  of $C(X)$, which are generated by two
subsets of idempotents, $Ann(I)+Ann(J)=Ann(I\cap J)$ (Question 1)
\ifif $X$ is an $EF$-space \ifif $\beta X$ is an $EF$-space.
%-----------------------------------------------------------
\begin{lem}
Let $X$ be a  normal space. The following statements are
equivalent.
\begin{itemize}
\item[$(a)$] $X$ is an $EF$-space. \item[$(b)$] For any two
collections $\mathcal{U}$ and $\mathcal{V}$  of clopen subsets of
$X$ with $\bigcup \mathcal{U}\cap \bigcup \mathcal{V}=\emptyset$,
we have $cl(\bigcup \mathcal{U})\cap cl(\bigcup
\mathcal{V})=\emptyset.$
\end{itemize}
\end{lem}
\textbf{Proof.} $(a)\Leftrightarrow(b)$ This follows from the
definition of $EF$-space and the fact that in a normal space
disjoint closed subsets are completely separated.$\hfill\square$
%-------------------------------------------------------------------
\begin{prop}\label{3.2}
For a topological space $X$ consider the following statements:
\begin{itemize}
\item[$(a)$] Any union of clopen subsets of $X$ is
$C^{*}$-embedded. \item[$(b)$] The closure of any union of clopen
subsets of $X$ is an open subset. \item[$(c)$] $X$ is an
$EF$-space. \item[$(d)$] Disjoint unions of clopen subsets have
disjoint closures\item[$(e)$] If $O_{1}$ and $O_{2}$ are disjoint
open subsets  with $O_{1}$ a union of clopen sets, then their
closures are disjoint.\\\indent Then
$(a)\Rightarrow(c)\Rightarrow(d)$ and
$(b)\Leftrightarrow(e)\Rightarrow(c)$.
\end{itemize}
\end{prop}
\textbf{Proof.} $(a)\Rightarrow(c)$ Let $\mathcal{U}$ and
$\mathcal{V}$ be two collections of clopen subsets and $\bigcup
\mathcal{U}\cap \bigcup \mathcal{V}=\emptyset$. Define $f(x)=1$
for $x\in \bigcup \mathcal{U}$ and $f(x)=0$ for $x\in \bigcup
\mathcal{V}$. Then $f\in C^{*}(\bigcup\mathcal{U}\cup \bigcup
\mathcal{V})$. By hypothesis, there exists $g\in C^{*}(X)$ such
that $g|_{(\bigcup\mathcal{U}\cup \bigcup \mathcal{V})}=f$.
Therefore $\bigcup \mathcal{U}\sub Z(1-g)$ and $\bigcup
\mathcal{V}\sub Z(g)$, i.e, they are completely
separated.\\\indent $(c)\Rightarrow(d)$ Trivial.
\\\indent $(b)\Leftrightarrow(e)$ This  follows from
\cite[Proposition 3.29]{KM}.
\\\indent $(e)\Rightarrow(c)$  Assume that $\mathcal{U}$ and
$\mathcal{V}$ are two collections of clopen subsets and $\bigcup
\mathcal{U}\cap \bigcup \mathcal{V}=\emptyset$.  Then
$cl(\bigcup\mathcal{U})\cap cl(\bigcup \mathcal{V})=\emptyset$.
Because $cl(\bigcup \mathcal{U})$ and $cl(\bigcup\mathcal{V})$ are
clopen subsets. Hence $cl(\bigcup \mathcal{U})$ and $cl(\bigcup
\mathcal{V})$ are disjoint zero sets, so $\bigcup \mathcal{U}$
 and $\bigcup
\mathcal{V}$ are completely separated. $\hfill\square$
%------------------------------------------------------------------
\\\indent In the following example, we see that $(a)$ does not imply $(b)$ in general.
This example was suggested to the author by A. Dow.
\begin{exam}
 Let   $\Bbb N$  denote the integers  and let $p$ be an ultra-filter on $\Bbb N$. Consider $X=\Bbb N\cup \{p\}\cup [0,1]$. Now
take the quotient space of $X$ in which the point  $p $  is
identified with $0$ in $[0,1]$. Then, we can see that any union of
clopen subsets is $C^{*}$-embedded. On the other hand $\Bbb N$ is
the union of clopen sets  but  its closure  includes $0$ from
$[0,1]$ and so is not clopen.
\end{exam}

\begin{prop}
If the closure of any union of clopen subsets is open and any
union of clopen subsets is $C^{*}$-embedded in its closure, then
any union of clopen subsets of $X$ is $C^{*}$-embedded in $X$.
\end{prop}
\textbf{Proof.} Let $\mathcal{U}$ be a collection of clopen
subsets and $A$ and $B$ be completely separated in $\bigcup
\mathcal{U}$. By Theorems \ref{1.2} and \ref{1.3}, it is enough to
prove that they are contained in disjoint zero-sets in $X$. By
hypothesis, there are $g_{1}, g_{2}\in C^{*}(cl(\bigcup
\mathcal{U}))$ such that $A\sub Z(g_{1})$, $B\sub Z(g_{2})$ and
$Z(g_{1})\cap Z(g_{2})=\emptyset$. Again by hypothesis, there
exists an idempotent element $f\in C^{*}(X)$ such that $cl(\bigcup
\mathcal{U}))=coz(f).$ Now we define $h_{1}(x)=0$, for all
$x\notin coz(f)$ and $h_{1}(x)=(g_{1}f)(x)$, for all $x\in
coz(f)$. Also, define $h_{2}(x)=1$, for all $x\notin coz(f)$ and
$h_{2}(x)=(g_{2}f)(x)$, for all $x\in coz(f)$. Then we have
$h_{1}, h_{2}\in C^{*}(X)$, $A\sub Z(g_{1})\sub Z(h_{1})$, $B\sub
Z(g_{2})\sub Z(h_{2})$ and $Z(h_{1})\cap Z(h_{2})=\emptyset$.
$\hfill\square$
%------------------------------------------------------------------------------
%-----------------------------------------------------------
\\\indent In \cite{C}, Comfort and Negrepontis restrict their attention to the
class of zero-dimensional space. For a space $X$ and open subset
$G$ of $X$, $G$ is of $type <\alpha$ if $G$ is the union of
$<\alpha$ closed-and-open subsets of $X$. Then, in the sense of
\cite{C}, a space $X$ is an $\textit{F}_{\alpha}$-$\textit{space}$
if every open subset of $X$ of $type <\alpha$ is $C^{*}$-embedded
in $X$ \cite[pp. 350, 343]{C}. So, by Proposition \ref{3.2} we
have the following result.
\begin{cor}\label{c}
If $X$ is an $F_{\alpha}$-space in the sense of \cite{C} for each
infinite cardinal $\alpha$, then $X$ is an $EF$-space.
\end{cor}
%---------------------------------------
\begin{exam}
A closed subset of an $EF$-space need not be an $EF$-space. To see
this, let $X=\beta \Bbb N\setminus\Bbb N$ as a closed subspace of
the $EF$-space $\beta \Bbb N$. Then by \cite[6S]{G}, the sets
$A^{\prime}=cl_{\beta\Bbb N}A\setminus \Bbb N$ form a base for the
open sets in $\beta \Bbb N\setminus\Bbb N$. So $X$ is a
zero-dimensional space. On the other hand by \cite[6W, 3]{G}, $X$
is not extremally disconnected. Hence, $X$ is not an $EF$-space.
\end{exam}
%------------------------------------------------------------%--------------------------------------------------------
The following example shows that  a $P$-space and hence a
basically disconnected space need not be an $EF$-space.
\begin{exam}\cite[4. N]{G}
Let $X$ be an uncountable space in which all points are isolated
points except for a distinguished point $s$. A neighborhood of $s$
is any set containing $s$ whose complement is countable, so any
set containing $s$ is closed. It is easily seen that $X$ is a
$P$-space. So $X$ is  basically disconnected. Now, consider two
disjoint uncountable subsets $U, V\sub X\setminus\{s\}$. Then we
have $s\in clU\cap clV$. So $X$ is not an $EF$-space.
\end{exam}
%-------------------------------------------------
%\begin{prop}
%Let $X=\bigoplus_{\alpha\in S}X_{\alpha}$ (see \cite[2.2]{E}).
%Then $X$ is an $EF$-space \ifif each  $X_{\alpha}$ is an
%$EF$-space.
%\end{prop}
%\textbf{Proof.} For each $\alpha\in S$, $X_{\alpha}$ is clopen in
%$X$. Hence, if $X$ is an $EF$-space, then each $X_{\alpha}$ is an
%$EF$-space. Now let $X_{\alpha}$ is an $EF$-space, for each
%$\alpha\in S$. $X_{\alpha}$ is $C^{*}$-embedded in $X$. So by
%Theorem \ref{1.2}, any two completely separated in $X_{\alpha}$
%are completely seperated in $X$, i.e, $X$ is an $EF$-space.
%-----------------------------------------------------
\begin{thm}\label{3.3}
A space $X$ is an $EF$-space \ifif $\beta X$ is an $EF$-space.
\end{thm}
\textbf{Proof.} Let $\{A_{\alpha}:\alpha\in S\}$ and
$\{A_{\gamma}:\gamma\in K\}$ be two collections of clopen subsets
of $\beta X$ and
\begin{center}
$(\bigcup_{\alpha\in S}A_{\alpha})\cap(\bigcup_{\gamma\in
K}A_{\gamma})=\emptyset$.
\end{center}
Then we have
\begin{center}
$(\bigcup_{\alpha\in S}A_{\alpha}\cap X)\cap(\bigcup_{\gamma\in
K}A_{\gamma}\cap X)=\emptyset.$
\end{center}
%\vspace{-.30cm}
 By hypothesis, there are
disjoint zero-sets $Z_{1}$ and $Z_{2}$ in $Z[X]$ such that;
%\vspace{-.30cm}
\begin{center}
%\begin{eqnarray*}
 $(\bigcup_{\alpha\in S}A_{\alpha})\cap X\sub Z_{1}$
and $(\bigcup_{\gamma\in K}A_{\gamma})\cap X\sub Z_{2}.$
\end{center}
%\end{eqnarray*}
 %\vspace{-.30cm}
Therefore we have,
%\vspace{-.30cm}
\begin{center}
%\begin{eqnarray*}
$\bigcup_{\alpha\in S}A_{\alpha}\sub cl_{\beta
X}(\bigcup_{\alpha\in S}A_{\alpha})=cl_{\beta
X}((\bigcup_{\alpha\in S}A_{\alpha})\cap X)\sub cl_{\beta
X}Z_{1}$,
%\end{eqnarray*}
\end{center}
%\vspace{-.55cm}
\begin{center}
 $\bigcup_{\gamma\in
K}A_{\gamma}\sub cl_{\beta X}(\bigcup_{\beta\in
K}A_{\gamma})=cl_{\beta X}((\bigcup_{\gamma\in K}A_{\gamma})\cap
X)\sub cl_{\beta X}Z_{2}$.
\end{center}
%\vspace{-.3cm}
 On the other hand $Z_{1}\cap Z_{2}=\emptyset$ implies that $cl_{\beta X}Z_{1}\cap cl _{\beta
X}Z_{2}=\emptyset$. By normality of $\beta X$, $cl_{\beta X}Z_{1}$
and $cl_{\beta X}Z_{2}$ are completely separated, i.e,
$\bigcup_{\alpha\in S}A_{\alpha}$ and $\bigcup_{\gamma\in
K}A_{\gamma}$ are completely separated.\\\indent Conversely,
assume that $\{A_{\alpha}:\alpha\in S\}$, $\{A_{\gamma}:\gamma\in
K\}$ are two collections of clopen subsets of $X$ and
\begin{center}
$(\bigcup_{\alpha\in S}A_{\alpha})\cap(\bigcup_{\gamma\in
K}A_{\gamma})=\emptyset$.
\end{center}
By Lemma \ref{1.1}, for
each $\alpha\in S$ and $\gamma\in K$, $cl_{\beta X}A_{\alpha}$ and
$cl_{\beta X}A_{\gamma}$ are clopen subsets of $\beta X$, so
\begin{center}
$(\bigcup_{\alpha\in S}cl_{\beta
X}A_{\alpha})\cap(\bigcup_{\gamma\in K}cl_{\beta
X}A_{\gamma})=\emptyset.$
\end{center}
 This together with our  hypothesis
implies that $\bigcup_{\alpha\in S}cl_{\beta X}A_{\alpha}$ and
$\bigcup_{\gamma\in K}cl_{\beta X}A_{\gamma}$ are completely
separated in $\beta X$. Hence $\bigcup_{\alpha\in S}A_{\alpha}$
and $\bigcup_{\gamma\in K}A_{\gamma}$ are completely separated in
$\ X$. $\hfill\square$
%-------------------------------------------------
%\\\indent We need the following lemma which is in $\cite[6R]{G}$.
%\begin{lem}\label{3.4}
%The following statements are equivalent.\\\indent $(a)$ Every
%subspace of $X$ is $C^{*}$-embedded.\\\indent $(b)$ $X$ is normal
%and every closed subspace is extremally disconnected.\\\indent
%$(c)$ $X$ is normal and every subspace is extremally disconnected.
%\end{lem}
%---------------------------------------------%------------------------------------------------------------
\\\indent In the next theorem we will answer Question 1. This result is an algebraic characterization of a completely
regular Hausdorff $EF$-space. First we need the following lemma.
\begin{lem}\label{3.5}
Let $R$ be a reduced ring. The following statements are
equivalent.
\begin{itemize}
\item[$(a)$] For any two orthogonal ideals $I$ and $J$ of $R$,
$Ann(I)+Ann(J)=R$.\item[$(b)$] For any two ideals $I$ and $J$ of
$R$, $Ann(I)+Ann(J)=Ann(I\cap J)$.
\end{itemize}
\end{lem}
\textbf{Proof.} $(a)\Rightarrow (b)$. We always have
$Ann(I)+Ann(J)\sub Ann(I\cap J)$. Now suppose that
 $x\in
Ann(I\cap J)=Ann(IJ)$. Then $xIJ=0$. So by  $(a)$ we have,
\begin{center}
$Ann(xI)+Ann(J)=R.$
\end{center}
 This shows that $1=a+b$, where $a\in Ann(xI)$
and $b\in Ann(J)$. Therefore $x=xa+xb$, where $xa\in Ann(I)$ and
$xb\in Ann(J)$ that is; $Ann(I\cap J)\sub Ann(I)+Ann(J).$\\\indent
$(b)\Rightarrow (a)$. If $IJ=0$, then by hypothesis,  we have
%\vspace{-.10cm}
\begin{center}
 $Ann(I)+Ann(J)=Ann(I\cap
J)=Ann(IJ)=R.$
\end{center}$\hfill\square$
%-------------------------------------------------------------%--------------------------------------------------------
\\\indent Note that, for any subset $S$ of $C(X)$,  $\bigcup COZ(S)=\bigcup_{f\in S}coz(f)$.
\begin{thm}\label{3.6}
Let $X$ be a completely regular Hausdorff space. The following
statements are equivalent.
\begin{itemize}
\item[$(a)$] The space $X$ is an $EF$-space.\item[$(b)$] For any
two ideals $I$ and $J$ of $C(X)$ which are generated by two
subsets of idempotents of $C(X)$, $Ann(I)+Ann(J)=Ann(I\cap J)$.
\item[$(c)$] For any two ideals $I$ and $J$ of $C^{*}(X)$ which
are generated by two subsets of idempotent elements of $C^{*}(X)$,
$Ann(I)+Ann(J)=Ann(I\cap J)$.
\end{itemize}
\end{thm}
\textbf{Proof.} $(a)\Rightarrow(b)$ Assume that $I$ and $J$ are
ideals of $C(X)$  generated by two subsets $S_{1}$ and $S_{2}$ of
idempotents. By Lemma \ref{3.5}, we can let $IJ=0$ and it is
enough to prove that $Ann(I)+Ann(J)=C(X)$. Now $IJ=0$ implies that
%\vspace{-.10cm}
\begin{center}
$(\bigcup COZ(I))\cap (\bigcup COZ(J))= (\bigcup COZ(S_{1}))\cap
(\bigcup COZ(S_{2}))=\emptyset.$
\end{center}
%\vspace{-.10cm}
 So there are disjoint zero-sets
$Z(f_{1})$ and $Z(f_{2})$ such that;
% \vspace{-.10cm}
\begin{center}
$\bigcup COZ(I)=\bigcup COZ(S_{1})\sub Z(f_{1})$ and $\bigcup
COZ(J)=\bigcup COZ(S_{2})\sub Z(f_{2})$.
\end{center}
%\vspace{-.10cm}
 Therefore $f_{1}\in
Ann(I)$, $f_{2}\in Ann(J)$ and $Z(f_{1}^{2}+f_{2}^{2})=\emptyset$.
Hence $f_{1}^{2}+f_{2}^{2}$ is a unit element in $Ann(I)+Ann(J)$,
i.e, $Ann(I)+Ann(J)=C(X)$. \\\indent $(b)\Rightarrow(a)$ Suppose
that $(\bigcup_{\alpha\in S}A_{\alpha})\cap(\bigcup_{\beta\in
K}A_{\beta})=\emptyset$, where for each $\alpha\in S$, $\beta\in
K$, $A_{\alpha}$ and $A_{\beta}$ are clopen subsets of $X$. It is
easily seen that for each $\alpha\in S$ and $\beta\in K$, there
are idempotent elements $e_{\alpha}$ and $e_{\beta}$ in $C(X)$
such that $A_{\alpha}=coz(e_{\alpha})$ and
$A_{\beta}=coz(e_{\beta})$. Hence
\begin{center}
$(\bigcup_{\alpha\in S}coz(e_{\alpha}))\cap(\bigcup_{\beta\in
K}coz(e_{\beta}))=\emptyset$.
\end{center}
 Now assume that $I=<\{e_{\alpha}:\alpha\in S\}>$ and $J=<\{e_{\beta}:\beta\in
 K\}>$. Then we can see that $IJ=0$. By hypothesis and Lemma \ref{3.5},
 $Ann(I)+Ann(J)=C(X).$
Hence there are $f\in Ann(I)$ and $g\in Ann(J)$ such that $1=f+g$.
Therefore we have,
%\vspace{-.1cm}
\begin{center}
 $\bigcup_{\alpha\in
S}A_{\alpha}=\bigcup_{\alpha\in S}coz(e_{\alpha})=\bigcup
COZ(I)\sub Z(f)$,
\end{center}
% \vspace{-55cm}
and
\begin{center}
 $\bigcup_{\beta\in K}A_{\beta}=\bigcup_{\beta\in
K}coz(e_{\beta})=\bigcup COZ(J)\sub Z(g).$
\end{center}
%\vspace{-.1cm}
 On the other hand $Z(f)\cap
Z(g)=\emptyset$. Hence $X$ is an $EF$-space.\\\indent
$(c)\Leftrightarrow(a)$ This is a consequence of
$(a)\Leftrightarrow(b)$, Theorem \ref{3.3}  and the fact that
$C^{*}(X)$ is isomorphic to $C(\beta X)$. $\hfill\square$
%---------------------------------------------------------
\\\indent Recall that, for a subset $A$ of $X$, we have $M_{A}=\{f\in C(X):
A\sub Z(f)\}$.
\begin{lem}\label{2.4}
\begin{itemize}
\item[$(a)$] For every subset $S$ of $C(X)$  we have,
$Ann(S)=M_{(\bigcup COZ[S])}$. \item[$(b)$] For subsets $A, B$ of
$X$ we have, $cl_{X}{A}=cl_{X}{B}$ \ifif $M_{A}=M_{B}$.
\end{itemize}
\end{lem}
\textbf{Proof.} $(a)$ Let $f\in Ann(S)$. Then $fg=0$, for all
$g\in S$. This implies that $\bigcup Coz[S]\sub Z(f)$, i.e, $f\in
M_{(\bigcup Coz[S])}$. Now $f\in M_{(\bigcup Coz[S])}$, implies
that, $Coz(g)\sub \bigcup Coz[S]\sub Z(f)$ for each $g\in S$, so
$f\in Ann(S)$.

$(b)$  If $cl_{X}{A}=cl_{X}{B}$, then it is easily seen that
$M_{A}=M_{cl_{X}A}=M_{cl_{X}B}=M_{B}$. Now suppose that
$M_{A}=M_{B}$ and $x\in cl_{X}{A}$. Then $M_{A}=M_{cl_{X}{A}}\sub
M_{x}$. If $x\notin cl_{X}{B}$, then by completely regularity of
$X$, there exists $f\in C(X)$ such that $x\notin Z(f)$ and $B\sub
Z(f)$, i.e, $f\in M_{B}\setminus M_{x}$,
 a contradiction. Hence $cl_{X}{A}\sub cl_{X}{B}$. Similarly,
we can prove that $cl_{X}{B}\sub cl_{X}A$. $\hfill\square$
%------------------------------------------------------------------
\\\indent From \cite{K}, a commutative ring $R$ with identity is  Baer if the
annihilator of every nonempty subset  of $R$ is generated by an
idempotent. The next result is proved in \cite{AZ}. Now we give a
new proof using Lemma \ref{2.4}.
\begin{thm}\label{3.57}\cite[Theorem 3.5]{AZ}
$C(X)$ is a Baer ring \ifif $X$ is an extremally disconnected
space.
\end{thm}
\textbf{Proof.} Let $A$ be an open subset of $X$. Then by
completely regularity of $X$, there exists a subset $S$ of $C(X)$
such that $A=\bigcup COZ[S]$. By hypothesis, there is an
idempotent $e\in C(X)$ such that $Ann(S)=Ann(e)$. By Lemma
\ref{2.4}, $M_{(\bigcup COZ[S])}=M_{coz(e)}$. Hence
\begin{center}
$clA=cl(\bigcup COZ[S])=cl(coz(e))=coz(e)$ is open.
\end{center}
 Conversely,
suppose that $S$ is a subset of $C(X)$. Then by hypothesis,
$cl(\bigcup COZ[S])$ is open. So there exists an idempotent $e\in
C(X)$ such that $cl(\bigcup COZ[S])=coz(e).$ Again by Lemma
\ref{2.4}, we have,
%\begin{center}
$Ann(S)=M_{(\bigcup COZ[S])}=M_{coz(e)}=Ann(e)$. $\hfill\square$
%\end{center}
%--------------------------------------------------------------------
\\\indent Recall that a commutative ring $R$ is an $
\textit{SA}$-$\textit{ring}$ (resp.,
$\textit{IN}$-$\textit{ring}$) if for any two ideals $I$ and $ J$
of $R$,
%\begin{center}
$Ann(I)+Ann(J)=Ann(K)$, for some ideal $K$ of $R$ (resp.,
$Ann(I)+Ann(J)=Ann(I\cap J)$), (see \cite{B}).
%\end{center}
\begin{cor}\label{3.8}
The following statements are equivalent.
\begin{itemize}
\item[$(a)$] $X$ is an extremally disconnected space.\item[$(b)$]
$C(X)$ is an $IN$-ring. \item[$(c)$] $C(X)$ is an
$SA$-ring.\item[$(d)$] The space of prime ideals of $C(X)$ is an
extremely disconnected space.
\end{itemize}
\end{cor}
\textbf{Proof.} This is a consequence of \cite[Corollary, 4.5]{B},
and Theorem \ref{3.57}. $\hfill\square$
%-----------------------------------------------------------------------------------------------------
\\\indent In \cite{M}, Swardson introduced an $\alpha$-open subset as a set of the form $A=\bigcup\mathcal{U}$, where $\mathcal{U}$ is a
collection of cozero-sets of $X$ with $|\mathcal{U} |< \alpha$.
She Also defined an $F_{\alpha}$-space $X$ to be a Tychonoff space
in which every $\alpha$-open subset of $X$ is $C^{*}$-embedded in
$X$. She proved that a space $X$ is an $F_{\alpha}$-space if and
only if any two disjoint $\alpha$-open subsets of $X$ are
completely separated in $X$ (see \cite[Theorem 2.3]{M}).
%---------------------------------------------------------------
\begin{thm}\label{3.9}
A Tychonoff space $X$ is an $F_{\alpha}$-space in the sense of
\cite{M} \ifif for any two $\alpha$-generated ideals $I$ and $J$
of $C(X)$, $Ann(I)+Ann(J)=Ann(I\cap J)$.
\end{thm}
\textbf{Proof.} By \cite[Theorem 2.3]{M}, the proof is similar to
that of Theorem \ref{3.6}.$\hfill\square$
%----------------------------------------------------------
Recall that, we mean of $\omega_{1}$-generated ideal is a
countably generated ideal. We also note that if $I=(f_{1},...,
f_{n},...)$ is a countably generated ideal in $C(X)$, then
$\rm{ann(I)}= \rm{ann(f)}$, where
$f=\sum_{i=1}^{\infty}\frac{|f_{n}|}{2^{n}(|f_{n}|+1)}.$
\begin{cor}
A topological space $X$ is an $F$-space \ifif for any two $f, g\in
C(X)$, $Ann(f)+Ann(g)=Ann(fg)$.
\end{cor}
\textbf{Proof.} ($\Rightarrow$) Let $fg=0$. Then $coz(f)\cap
coz(g)=\emptyset$. By hypothesis, there are $f_{1}, f_{2}\in C(X)$
such that $coz(f)\sub z(f_{1})$, $coz(g)\sub z(f_{2})$ and
$z(f_{1})\cap z(f_{2})=\emptyset$ so $f_{1}\in Ann(f)$, $f_{2}\in
Ann(g)$ and $f_{1}^{2}+f_{2}^{2}$ is a unit element in
$Ann(f)+Ann(g)$. Thus $Ann(f)+Ann(g)=C(X)$. By this equality and
Lemma \ref{3.5}, we have
\begin{center}
$Ann(f)+Ann(g)=Ann(<f>)+Ann(<g>)=Ann(<f>\cap <g>)=Ann(fg)$.
\end{center}
($\Leftarrow$) Clearly $X$ is an $F_{\omega_{1}}$-space in the
sense of \cite{M} if and only if $X$ is an F-space. Now let $I$ be
a ${\omega_{1}}$-generated ideal. By the above comment, there
exists $f\in C(X)$ such that $Ann(I)=Ann(f)$. By Theorem
\ref{3.9}, $X$ is an $F$-space.$\hfill\square$
%--------------------------------------------------------
\begin{cor}\cite[Proposition,
2.2]{M}\label{3.10} A completely regular Hausdorff space $X$ is an
extremely disconnected space \ifif $X$ is an $F_{\alpha}$-space in
the sense of \cite{M} for each infinite cardinal $\alpha$ .
\end{cor}
\textbf{Proof.} This is a consequence of Corollary \ref{3.8} and
Theorem \ref{3.9}.
%---------------------------------------------------------------

%\begin{lem}
%If $A$ is an open subset and $B$ is any subset of $X$, then
%$intcl_{X}({A\cap B})=intcl_{X}{A}\cap intcl_{X}{B}$.
%\end{lem}
%--------------------------------------------------------
%\begin{lem}
%A subset $F$ of $Spec(R)$ is an irreducible closed   \ifif
%$F=V(P)$, for some prime ideal $P$ of $R$.
%\end{lem}
%\textbf{Proof}. Let $F$ be closed subset of $Spec(R)$. Then there
%is an ideal $I$ of $R$ such that $F=V(I)=V(\sqrt{I})$. Now let $a,
%b\in R$ and $ab\in \sqrt{I}$. Then $F=V(I)\sub V(a)\cup V(b)$, so
%$F=V(I)\sub V(a)$ or $F=V(I)\sub V(b)$, i.e, $a\in \sqrt{I}$ or
%$b\in\sqrt{I}$. Therefore $F=V(\sqrt{I})$, where $\sqrt{I}$ is
%prime. Conversely, Let $F\sub F_{1}\cup F_{2}$, for some closed
%subsets $F_{1}$ and $F_{2}$ of $Spec(R)$. Then there are ideals
%$I$ and $J$ of $R$ such that  $F=V(P)\sub V(I)\cup V(J)=V(IJ)$.
%This implies that $IJ\sub\sqrt{IJ}\sub P$, hence $I\sub P$ or
%$J\sub P$, i.e, $F=V(P)\sub V(I)=F_{1}$ or $F=V(P)\sub
%V(J)=F_{2}$.
%---------------------------------------------------

\section{$EZ$-spaces}
%These are the main results of the paper.
%\begin{defn}
%A topological space $X$ with a base for open subsets which the
%closure of each  element is open denote by $E_{1}$-space
%\end{defn}
%------------------------------------------------
%\begin{defn}
%A topological space $X$ which has a base $\cal B$ for open subsets
%such that $U, V\in\cal B$ and $U\cap V=\emptyset$, implies that
%$\overline{U}\cap \overline{V}=\emptyset$ denote by $F_{1}$-space.
%\end{defn}

%-------------------------------------------------------------
We call a topological space $X$ an $\textit{EZ}$-$\textit{space}$
if for every open subset $A$ of $X$ there exists a collection
$\{A_{\alpha}: \alpha\in S\}$ of clopen subsets of $X$ such that
$cl_{X}{A}=cl_{X}({\bigcup_{\alpha\in S}A_{\alpha}})$,
equivalently, every regular closed subset is the closure of a
union of clopen sets. In this section, it is proved that a
completely regular Hausdorff space $X$ is an $EZ$-space \ifif for
any ideal $I$ of $C(X)$ there exists a subset $S$ of idempotents
of $C(X)$ such that $Ann(I)=Ann(S)$ (Question 2) \ifif for any
$f\in C(X)$ there exists a subset $S$ of idempotents of $C(X)$
such that $Ann(f)=Ann(S)$ \ifif $\beta X$ is an $EZ$-space. As a
consequence, a completely regular Hausdorff space $X$ is an $EF$
and $EZ$-space \ifif $X$ is an extremally disconnected space \ifif
$X$ is an $F_{\alpha}$-space for each infinite cardinal $\alpha$
in the sense of \cite{M} \ifif $X$ is an $F_{\alpha}$-space for
each infinite cardinal $\alpha$ in the sense of \cite{C}.
%----------------------------------------------------------
\\\indent Recall that a collection $\cal{B}$ of open sets in a topological space
$X$ is called a $\pi$-base if every open set of $X$ contains a
member of $\cal{B}$. Thus, a clopen $\pi$-base is one consisting
of clopen sets. The reader is referred to \cite{BH}, \cite{BM},
\cite{E}, \cite{KM}, \cite{MC} and \cite{V}.
 In Proposition \ref{z}, we find equivalent conditions for when a space with a clopen
 $\pi$-base is an $EF$-space.
\begin{prop}\label{2.1}
\begin{itemize}
\item[$(a)$] Every extremally disconnected space is an
$EZ$-space.\item[$(b)$] If $X$ has a clopen $\pi$-base, then $X$
is an $EZ$-space.\item[\rm(c)] Every $T_{1}$-space with  a dense
set of isolated points is an $EZ$-space.\item[$(d)$] Every open
subset of an $EZ$-space is an $EZ$-space.
\end{itemize}
\end{prop}
\textbf{Proof.} $(a)$ is  obvious.\\\indent $(b)$ Let $B$ be an
open subset of $X$. Set $\mathcal{O}$ as the collection of clopen
subsets of $X$ contained in $cl_{X}B$. We claim that
$cl_{X}(\bigcup\mathcal{O})=cl_{X}B$. Clearly
$cl_{X}(\bigcup\mathcal{O})\sub cl_{X}B$. Suppose there is an
$x\in cl_{X}B\setminus cl(\bigcup\mathcal{O})$. Since $cl_{X}B$ is
regular closed set $x\in cl_{X}int_{X}cl_{X}B$, so that
$int_{X}cl_{X}B\cap (X\setminus cl_{X}(\bigcup\mathcal{O}))$ is
nonempty. By hypothesis, there is a nonempty clopen subset in
$int_{X}cl_{X}B\cap (X\setminus cl_{X}(\bigcup\mathcal{O}))$.
However, such a clopen set is contained in $cl_{X}B$ and
$(X\setminus cl_{X}(\bigcup\mathcal{O}))$, a contradiction.
\\\indent $(c)$ A $T_{1}$-space with a dense
set of isolated points has a clopen $\pi$-base and hence is an
$EZ$-space.\\\indent (d) Assume that $Y$ is an open subset of $X$.
Then by hypothesis, for any open subset $A$ of $Y$, there exists a
collection $\{A_{\alpha}: \alpha\in S\}$ of clopen subsets of $X$
such that $cl_{X}(A)=cl_{X}(\bigcup_{\alpha\in S}A_{\alpha})$. Now
$Y$ is open in $X$, so
\begin{center}
$cl_{Y}(A)=cl_{X}(A)\cap Y=cl_{X}(\bigcup_{\alpha\in
S}A_{\alpha})\cap Y=cl_{Y}(\bigcup_{\alpha\in S}A_{\alpha}\cap
Y)$.
\end{center}$\hfill\square$
%----------------------------------------------------------------
\begin{lem}\label{ref1}
Let $X$ be a regular space. Then $X$ is an $EZ$-space \ifif $X$
has a clopen $\pi$-base.
\end{lem}
\textbf{Proof.} By Proposition \ref{2.1}, the necessity is
evident. Suppose $X$ is an $EZ$-space and let $O$ be any nonempty
open subset of $X$. Choose $x\in O$ and by regularity choose an
open neighborhood of $x$, say $T$, such that $x\in T\sub clT\sub
O$. By hypothesis, $clT = clU$ where $U$ is a union of clopen
subsets. Then any nonempty clopen subset of $U$ is also a nonempty
clopen subset of $O$. Therefore, $X$ has a clopen $\pi$-base.
$\hfill\square$
%----------------------------------
\\\indent By Lemma \ref{ref1},  \cite[Proposition 18]{MC} and \cite[Corollary 19]{MC}, we have $X$ has clopen
$\pi$-base \ifif $\beta X$ has a clopen $\pi$-base.
\\\indent The following example shows that the regularity hypothesis,
in Lemma \ref{ref1} is not superfluous and hence an $EZ$-space
need not have a clopen $\pi$-base.
\begin{exam}\label{est}
Let $X=\Bbb R$ with the topology $\cal{T}$=$\{(a, \infty):
a\in\Bbb R\}\cup\{\emptyset, \Bbb R\}$. Then $X$ is a non-regular
space and any nonempty open subset is dense, so $X$ is an
$EZ$-space. On the other hand $X$ does not a clopen $\pi$-base.
\end{exam}
%--------------------------------------------------------------------%-------------------------------------------%-----------------------------------------------------
Example \ref{est} motivates the following question: Does there
exist a Hausdorff $EZ$-space with no clopen $\pi$-base?
%----------------------------------------------------------------
\\\indent Recall that from \cite{BM}, a $DC$-space is a Tychonoff space $X$
such that for each $f\in C(X)$ there exists a family of open
subsets $\{U_{i}:i\in I\}$, the union of which is dense in $X$,
such that $f$, restricted to each $U_{i}$ is constant. By
\cite[Lemma 2.5]{BH}, and Proposition \ref{2.1} any separable
$DC$-space is an $EZ$-space.
\begin{prop}\label{2.2}
 Let $X=\bigoplus_{\alpha\in S}X_{\alpha}$. Then $X$ is an $EZ$-space
 \ifif each
 $X_{\alpha}$ is an $EZ$-space.
 \end{prop}
 \textbf{Proof.} ($\Rightarrow$) Each $X_{\alpha}$ is open in $X$.
 By Proposition \ref{2.1} $(d)$, each $X_{\alpha}$ is an
 $EZ$-space.\\\indent
  ($\Leftarrow$) Suppose that $A$ is an open subset in $X$. Then $A\cap X_{\alpha}$ is open in $X_{\alpha}$. Therefore, for each $\alpha\in S$, there exists a
 collection $\{A^{\beta}_{\alpha}: \beta\in S_{\alpha}\}$ of clopen subsets
of $X_{\alpha}$ such that $cl_{X_{\alpha}}(A\cap
X_{\alpha})=cl_{X_{\alpha}}(\bigcup_{\beta\in
S_{\alpha}}A^{\beta}_{\alpha})$. Now it is easy to see  that
$cl_{X}A=cl_{X}\bigcup_{\alpha\in S}(\bigcup_{\beta\in
S_{\alpha}}A^{\beta}_{\alpha})$, i.e, $X$ is an
$EZ$-space.$\hfill\square$
%----------------------------------------------------------------
\begin{prop}\label{2.3}
%\begin{itemize}
 Every dense subset of an $EZ$-space is an $EZ$-space.
%\end{itemize}
\end{prop}
\textbf{Proof.} Let $Y$ be a dense subset of a topological space
$X$ and $A$ be an open subset of  $Y$. Then there exists an open
subset $G$ of $X$ such that $A=G\cap Y$. By hypothesis, there
exists a collection $\{A_{\alpha}: \alpha\in S\}$ of clopen
subsets of $X$ such that $cl_{X}{G}=cl_{X}({\bigcup_{\alpha\in
S}A_{\alpha}})$. On the other hand $cl_{X}A=cl_{X}(G\cap
Y)=cl_{X}G$. So we have,
%\begin{center}
$cl_{Y}A=cl_{X}A\cap Y=cl_{X}({\bigcup_{\alpha\in
S}A_{\alpha}}\cap Y)\cap Y =cl_{Y}(\bigcup_{\alpha\in
S}A_{\alpha}\cap Y)$. $\hfill\square$
%\end{center}
%-------------------------------------------------------
\\\indent A completely regular Hausdorff $EZ$-space need not even be
a totally disconnected space. The following example, (i.e., (a))
was presented by  Brian Scott for another purpose. Also, by
Proposition \ref{2.3}, we have another example of an $EZ$-space
which is not a zero-dimensional space.
\begin{exam}
\begin{itemize}
\item[\rm(a)] For $n\in\Bbb Z^+$ let
$$D_n=\{(\frac{2k+1}{2^n},\frac1{2^n}):k=0,\dots,2^{n-1}-1\},$$
and let $D=\bigcup_{n\in\Bbb Z^+}D_n\;.$ Now suppose that $X=D\cup
\{(a, 0): a\in [0,1]\}$ as a subspace of $\Bbb R^2$ with the usual
topology. Then $D$ is a countable dense subset of isolated points
in $X$. By Proposition \ref {2.1} $(c)$, $X$ is an $EZ$-space. On
the other hand $X$ is not totally disconnected  and hence X is not
zero-dimensional.\item[\rm(b)] By Dowker's Example \cite[Example,
6.2.20]{E}, we have a zero-dimensional space $Y$ for which $\beta
Y$ is not zero-dimensional. On the other hand, by Proposition
\ref{2.3}, $\beta Y$ is an $EZ$-space.
\end{itemize}
\end{exam}
A closed subset of an $EZ$-space need not be an $EZ$-space.
Because $Y=[0,1]$ as a closed subset of $X$ in the above example
(i.e., (a)),  is not an $EZ$-space.
%------------------------------------------------------------------

In the following theorem we answer Question $2$; it is an
algebraic characterization of a completely regular Hausdorff
$EZ$-space (i.e., a space with a clopen $\pi$-base).
%--------------------------------------------------------\Box\blacksquare
\begin{thm}\label{2.5} Let $X$ be a completely regular
Hausdorff space. The following statements are equivalent.
\begin{itemize}
\item[$(a)$] $X$ is an $EZ$-space (i.e., $X$ has a clopen
$\pi$-base). \item[$(b)$]  For every ideal $I$ of $C(X)$, there
exists a subset $S$ of idempotent elements of $C(X)$ such that
$Ann(I)=Ann(S)$. \item[$(c)$] For every $f\in C(X)$, there exists
a subset $S$ of idempotent elements of $C(X)$ such that
$Ann(f)=Ann(S)$. \item[$(d)$] For every cozero-set $H$ of $X$
there exists a collection $\{H_{\alpha}: \alpha\in S\}$ of clopen
subsets of $X$ such that $cl_{X}{H}=cl_{X}({\bigcup_{\alpha\in
S}H_{\alpha}})$.
\end{itemize}
\end{thm}
\textbf{Proof.} $(a)\Rightarrow(b)$ For an ideal $I$ of $C(X)$
consider, $A=\bigcup COZ[I]$. Then by hypothesis, there exists a
collection $\{A_{\alpha}: \alpha\in S\}$ of clopen subsets of $X$
such that $cl_{X}{A}=cl_{X}({\bigcup_{\alpha\in S}A_{\alpha}})$.
It is easily seen that for each $\alpha\in S$, there exists an
idempotent $e_{\alpha}$ such that $A_{\alpha}=coz(e_{\alpha})$.
Now suppose that $S=\{e_{\alpha}: \alpha\in S\}$. Then by Lemma
\ref{2.4}, we have,
\begin{center}
$Ann(I)=M_{\bigcup COZ[I]}=M_{\bigcup COZ[S]}=Ann(S).$
\end{center}
$(b)\Rightarrow(a)$ Let $A$ be an open subset of $X$. We know that
in a completely regular space $X$, $COZ[X]$ is a base for open
subsets. So there exists a subset $K$ of $C(X)$ such that
$A=\bigcup COZ[K]$. Now suppose that $I$ is the ideal generated by
$K$ in $C(X)$. Then by hypothesis, there exists a subset $S$ of
idempotent elements of $C(X)$ such that
$M_{A}=Ann(I)=Ann(S)=M_{\bigcup COZ[S]}$. Therefore by Lemma
\ref{2.4}, we have $cl_{X}(A)=cl_{X}(\bigcup COZ[S])$.\\\indent
$(b)\Rightarrow(c)$ For any $f\in C(X)$, we have
$Ann(f)=Ann(<f>)$. By hypothesis, there exists a subset $S$ of
idempotents such that $Ann(<f>)=Ann(S)$. Hence
$Ann(f)=Ann(S)$.\\\indent $(c)\Rightarrow(b)$ Let $I$ be an ideal
of $C(X)$. Then $Ann(I)=\bigcap_{f\in I}Ann(f)$. By hypothesis,
for each $f\in I$ there exists a subset $S_{f}$ of idempotent such
that $Ann(f)=Ann(S_{f})$. Therefore
\begin{center}
$Ann(I)=\bigcap_{f\in I}Ann(S_{f})=Ann(\bigcup_{f\in I}S_{f})$.
\end{center}
 $(c)\Leftrightarrow(d)$ This is similar to that of
$(a)\Leftrightarrow(b)$ step by step.$\hfill\square$
%----------------------------------------------
\\\indent In \cite{KM}, M. Knox and W. Wm. McGovern  define $X$ to be a $qsz$-space
if for any $f\in C(X)^{+}$ there exists a countable sequence
$K_{n}$ of clopen subsets $X$ such that
$cl_{X}coz(f)=cl_{X}\bigcup_{n\in\Bbb N}K_{n}$. So, by the above
theorem, we have  the following corollary.
%----------------------------------------------------------
\begin{cor}\cite[Proposition 3.9]{KM}
If $X$ is a $qsz$-space, then $X$ is an $EZ$-space.
\end{cor}
%----------------------------------------------------
 Recall that, a subspace $Y$ of a space $X$ is  $z$-embedded in
$X$ if for every zero-set $Z$ in $Y$ there is a zero-set $H$ in
$X$ such that $Z = H\cap Y$, equivalently,  for every cozeroset of
$Y$ there is a cozeroset of $X$ which traces to it. For example, a
$C^{*}$-embedded subspace is clearly $z$-embedded (see \cite{MA}).
\begin{cor}
If $X$ is a  completely regular Hausdorff $EZ$-space, then every
$z$-embedded subspace  is an $EZ$-space.
\end{cor}
\textbf{Proof.}  Let $Y$ be a $z$-embedded subspace of an
$EZ$-space $X$ and $H$ a cozero-set in $Y$. By hypothesis and
Theorem \ref{2.5}, there exists a cozero-set $C$ in $X$ and a
collection $\{A_{\alpha}:\alpha\in S\}$ of clopen subsets of $X$
such that $H=C\cap Y$ and $cl_{X}C=cl_{X}(\bigcup_{\alpha\in
S}A_{\alpha})$. Therefore
\begin{center}
$cl_{Y}H=cl_{Y}(C\cap Y)=cl_{X}C\cap Y=cl_{Y}(\bigcup_{\alpha\in
S}A_{\alpha}\cap Y)$.
\end{center}$\hfill\square$
%-----------------------------------------------
\begin{prop}\label{z}
Let $X$ be an $EZ$-space. Then the following statements are
equivalent.
\begin{itemize}
\item[$(a)$] $X$ is an $EF$-space. \item[$(b)$] If $\mathcal{U}$
and $\mathcal{V}$ are two collections of clopen subsets of $X$
with $\bigcup \mathcal{U}\cap \bigcup \mathcal{V}=\emptyset$, then
$cl(\bigcup \mathcal{U})\cap cl(\bigcup \mathcal{V})=\emptyset$.
 \item[$(c)$]
 The closure of
any union of clopen subsets of $X$ is an open subset. \item[$(d)$]
$X$ is an extremally disconnected space. \item[$(e)$] Any union of
clopen subsets of $X$ is $C^{*}$-embedded.
\end{itemize}
\end{prop}
\textbf{Proof.} By Proposition \ref{3.2} and definitions, it is
clear that $(d)\Rightarrow(e)\Rightarrow (a)\Rightarrow
(b)$.\\\indent $(b)\Rightarrow (c)$ Let $\cal{U}$ be a collection
of clopen subsets and let $B$ be an open subset of $X$ which is
disjoint from any element of $\cal{U}$. By \cite[Proposition
3.29]{KM}, it is enough to prove that $clB\cap
cl(\bigcup\cal{U})=\emptyset$. By hypothesis, there exists a
collection $\mathcal{O}$ of clopen subsets such that
$clB=cl(\bigcup\mathcal{O})$. We have $B\cap
(\bigcup\cal{U})=\emptyset$, so
\begin{center}
 $clB\cap
(\bigcup\mathcal{U})=cl(\bigcup\mathcal{O})\cap
\bigcup\mathcal{U}=\emptyset$.
\end{center}
 This shows that $\bigcup\mathcal{O}\cap
\bigcup\mathcal{U}=\emptyset$. Now  hypothesis implies that,
\begin{center}
$clB\cap cl(\bigcup\mathcal{U})=cl(\bigcup\mathcal{O})\cap
cl(\bigcup\mathcal{U})=\emptyset$.
\end{center}
 $(c)\Rightarrow(d)$ By
hypothesis, for any open subset $B$ there exists a collection
$\mathcal{O}$ of clopen subsets such that
$clB=cl(\bigcup\mathcal{O})$. Now, by hypothesis, $clB$ is open,
i.e, $X$ is extremally disconnected. $\hfill\square$
%----------------------------------------------
\\\indent It is well known that if $X$ is  zero-dimensional and an $F_{\alpha}$-space for each infinite cardinal $\alpha$ in the
sense of \cite{M}, then $X$ is an $F_{\alpha}$-space for each
infinite cardinal $\alpha$, in the sense of \cite{C}. In the next
result we see that the two concepts of $F_{\alpha}$-space coincide
for each infinite cardinal $\alpha$
%------------------------------------------------------------------------------------
\begin{cor}\label{2.7}
Let $X$ be a completely regular Hausdorff space. The following
statements are equivalent.
\begin{itemize}
\item[$(a)$] $X$ is an extremally disconnected space.\item[$(b)$]
$X$ is $EF$ and $EZ$-space. \item[$(c)$] $X$ is an
$F_{\alpha}$-space in the sense of \cite{M} for each infinite
cardinal $\alpha$.\item[$(d)$] $X$ is an $F_{\alpha}$-space in the
sense of \cite{C} for each infinite cardinal $\alpha$
.\item[$(e)$] The closure of any union of clopen subsets is open
and $X$ is an $EZ$-space.
\end{itemize}
\end{cor}
\textbf{Proof.} $(a)\Leftrightarrow(c)$ Follows from Corollary
\ref{3.10}.
\\\indent
$(a)\Rightarrow (b)$ By Propositions \ref{3.2} and \ref{2.1}, any
extremely disconnected space is  $EF$ and $EZ$.
\\\indent $(b)\Rightarrow (a)$ Let $X$ be  $EF$ and $EZ$-space.
Then by Proposition \ref{z}, $X$ is an extremally disconnected
space.
\\\indent $(a)\Rightarrow (d)$ Every extremally disconnected
$T_{3}$-space is zero-dimensional. On the other hand any open
subset of an extremally disconnected space is $C^{*}$-embedded so
any union of clopen subsets is $C^{*}$-embedded. Therefore $X$ is
an $F_{\alpha}$-space in the sense of \cite{C}  for each infinite
cardinal $\alpha$.\\\indent $(d)\Rightarrow (a)$ If $X$ is an
$F_{\alpha}$-space in the sense of \cite{C} for each infinite
cardinal $\alpha$, then by Corollary \ref{c}, $X$ is an
$EF$-space. On the other hand $X$ is zero-dimensional so is an
$EZ$-space. Hence by Proposition \ref{z}, $X$ is an extremally
disconnected space.
\\\indent $(a)\Rightarrow(e)$ This is obvious.\\\indent
$(e)\Rightarrow(a)$ Let $A$ be an open subset. Then by hypothesis,
there exists a collection $\{A_{\alpha}: \alpha\in S\}$ of clopen
subsets such that $clA=cl(\bigcup_{\alpha\in S}A_{\alpha})$. Again
by hypothesis, $clA$ is open, i.e, $X$ is an extremally
disconnected space.$\hfill\square$
%--------%-----------------------------------------------------------------

\section{$Spec(R)$ as an $EZ$-space}

In this section, for a reduced ring $R$, we prove that $Spec(R)$
is an $EZ$-space \ifif for every ideal $I$ of $R$ there exists a
subset $S$ of idempotents of $R$ such that $Ann(I)=Ann(S)$ (a
general case of Question 2) \ifif for any $a\in R$, there exists a
subset $S$ of idempotent elements of $R$ such that
$Ann(a)=Ann(S)$. Also, for a ring $R$ satisfying $J(R)=0$, we show
that $Max(R)$ is an $EZ$-space \ifif for every ideal $I$ of $R$
there exists a subset $S$ of idempotents of $R$ such that
$Ann(I)=Ann(S)$.
\\\indent For $a\in R$, let $supp(a)=\{P\in Spec(R): a\notin P\}$.
It is easy to see that for any $R$, $\{supp(a): a\in R\}$ forms a
basis of open sets for $Spec(R)$ (i.e., the space of prime ideals
of $R$). This topology is called the Zariski topology. We use
$V(I) (V(a))$ to denote the set of $P\in Spec(R)$, such that
$I\sub P (a\in P)$. Note that $V(I)=\bigcap_{a\in I}V(a)$ and
$V(a)=Spec(R)\setminus supp(a)$ (see \cite{HU}).
 \\\indent For an open subset $A$ of $Spec(R)$, let $O_{A}:=\{a\in R:
A\sub V(a)\}$. Since for any $a, b\in R$, $V(a)\cap V(b)\sub
V(a-b)$ and for each $r\in R$, $a\in O_{A}$, we have $V(a)\sub
V(ra)$,  thus $O_{A}$ is an ideal of $R$. It is easy to see that
$O_{A}=\bigcap_{P\in A}P$ and $V(O_{A})=clA$, where $clA$ is the
cluster points of $A$ in $Spec(R)$. \\\indent An ideal $I$ of a
commutative ring $R$ is said to be an annihilator ideal provided
that $Ann(Ann(I))=I$,  equivalently, if $Ann(I)\sub Ann(x)$, and
$x\in R$, then $x\in I$.
%-------------------------------------------
\\\indent We need the following lemma which consists of some well-known results.
\begin{lem}\label{2.7} Let $R$ be a reduced ring.
\begin{itemize}
\item[$(a)$] For ideals $I, J$ of $R$, $Ann(I)\sub Ann(J)$ \ifif
$intV(I)\sub intV(J).$ \item[$(b)$] For an open subset $A$ of
$Spec(R)$, $O_{A}$ is an annihilator ideal. \item[$(c)$] If $I$ is
an annihilator ideal of $R$, then there exists an open subset $A$
of $Spec(R)$ such that $I=O_{A}$. \item[$(d)$]For open subsets $A,
B$ of $Spec(R)$, $O_{A}= O_{B}$ \ifif $clA=clB$. \item[$(e)$]
$A\sub spec(R)$ is a clopen subset \ifif there is an idempotent
$e\in R$ such that $A=V(e)=supp(1-e)$. \item[$(f)$] For any ideal
$I$ of $R$, $Ann(I)=O_{(\bigcup supp(I))}$.
\end{itemize}
\end{lem}
\textbf{Proof.} $(a)$ let $I$ and $J$ be two ideals of $R$ and
$P\in intV(I).$ Then there is an $a\in R$ such that $P\in
supp(a)\sub V(I)$. Hence $supp(Ia)=supp(I)\cap supp(a)=\emptyset$,
thus $Ia=0$. This implies that $a\in Ann(I)\sub Ann(J)$, so
$Ja=0$. Therefore $P\in supp(a)\sub intV(J)$. Conversely, let
$x\in Ann(I)$. Then $Ix=0$. so $supp(x)\sub intV(I)\sub intV(J)$.
This shows that $supp(Jx)=supp(x)\cap supp(J)=\emptyset$. Hence
$Jx=0$, i.e, $x\in Ann(J)$.\\\indent $(b)$ Let $Ann(O_{A})\sub
Ann(x)$. By $(a)$,
\begin{center}
 $A\sub intclA=intV(O_{A})\sub intV(x)\sub V(x)$.
\end{center}
So $x\in O_{A}$.\\\indent
 (c) Suppose that $I$ is an annihilator ideal and
$A=intV(I)$. We claim that $I=O_{A}$. If $a\in I$, then
$intV(I)\sub V(a)$, i.e, $a\in O_{A}$. Now let $a\in O_{A}$. Then
$A=intV(I)\sub intV(a)$, so $Ann(I)\sub Ann(a)$. On the other hand
$I$ is an annihilator,  hence $a\in I$.
\\\indent $(d)$ We have $O_{A}= O_{B}$. This implies that $V(O_{B})= V(O_{A})$,
i.e., $cl(B)=cl(A)$. Conversely, $clA=clB$ implies that
$O_{A}=O_{clA}=O_{clB}=O_{B}$.\\\indent $(e)$  Let $A$ be a clopen
subset, $I=O_{A}$ and $J=O_{A^{c}}$. Then $A=clA=V(O_{A})=V(I)$
and $A^{c}=V(O_{A^{c}})=V(J)$. Hence $V(I+J)=V(I)\cap
V(J)=\emptyset$, so there are $a\in I$ and $b\in J$ such that
$1=a+b$. But $V(a)\cup V(b)=Spec(R)$, thus we have $ab=0$, this
implies that $a=a^{2}$ and $V(I)=V(a).$ The converse is
evident.\\\indent $(f)$ If $r\in Ann(I)$, then $ra=0$, for all
$a\in I$, so $\bigcup supp(I)\sub V(r)$, this shows that $r\in
O_{(\bigcup supp(I))}$. Now $r\in O_{(\bigcup supp(I))}$, implies
that $supp(a)\sub V(r)$, for all $a\in I$, so $supp(a)\cap
supp(r)=supp(ra)=\emptyset$, i.e, $ra=0$. Hence $r\in Ann(I)$.
$\hfill\square$
%---------------------------------------------------------
\begin{thm}\label{2.8}
Let $R$ be a reduced ring. The following statements are
equivalent.
\begin{itemize}
\item[$(a)$] The space of prime ideals, $Spec(R)$, is an
$EZ$-space.\item[$(b)$] For every ideal $I$ of $R$, there exists a
subset $E$ of idempotents of $R$ such that
$Ann(I)=Ann(E)$.\item[$(c)$] For every $a\in R$, there exists a
subset $S$ of idempotents of $R$ such that
$Ann(a)=Ann(S)$.\item[$(d)$] For any $a\in R$, there exists a
clopen subset $S$ of  $Spec(R)$ such that $cl(supp(a))=cl(\bigcup
supp(S))$.
\end{itemize}
\end{thm}
\textbf{Proof.} $(a)\Rightarrow(b)$ Let $I$ be an ideal of $R$.
Then we have $\bigcup supp(I)$ is an open subset of $Spec(R)$. By
hypothesis, there exists a collection $\{A_{\alpha}: \alpha\in
S\}$ of clopen subsets such that $cl(\bigcup
supp(I))=cl(\bigcup_{\alpha\in S} A_{\alpha})$. By Lemma
\ref{2.7}, for each $\alpha\in S$ there exists an idempotent
$e_{\alpha}$ such that $A_{\alpha}=supp(e_{\alpha})$. Therefore,
\begin{center}
$cl(\bigcup supp(I))=cl(\bigcup_{\alpha\in S} supp(e_{\alpha}))$.
\end{center}
Again by Lemma \ref{2.7}, we have $Ann(I)=Ann(E)$ where
$E=\{e_{\alpha}:\alpha\in S\}$.\\\indent $(b)\Rightarrow(a)$ Let
$A$ be an open subset of $Spec(R)$. Then there exists a subset $K$
of $R$ such that $A=\bigcup supp[K]$. Now suppose that $I$ be the
ideal generated by $K$ in $R$. Then by hypothesis and Lemma
\ref{2.7}, there exists a subset $E$ of idempotents of $R$ such
that $O_{A}=Ann(I)=Ann(E)=O_{\bigcup supp[E]}$. Therefore,  Lemma
\ref{2.7} implies that $cl(A)=cl(\bigcup supp[E])$.\\\indent
$(b)\Rightarrow(c)$ This is evident.\\\indent $(c)\Rightarrow(b)$
For an ideal $I$ of $R$ we have $Ann(I)=\bigcap_{a\in I}Ann(a)$.
By hypothesis, for each $a\in R$ there exists a subset $S_{a}$ of
idempotents such that $Ann(a)=Ann(S_{a})$. Hence
$Ann(I)=\bigcap_{a\in I}Ann(S_{a})=Ann(\bigcup_{a\in
I}S_{a})$.\\\indent $(c)\Leftrightarrow(d)$ By Lemma \ref{2.7},
$Ann(a)=Ann(S)$ for some subset $S$ of $R$ \ifif
$O_{supp(a)}=O_{\bigcup supp(S)}$ \ifif $cl(supp(a))=cl(\bigcup
supp(S))$. $\hfill\square$
%-----------------------------------------------
%-----------------------------------
\\\indent We denote by $Max(R)$ the space of maximal ideals of $R$. For $a\in R$, let $D(a)=\{M\in Max(R): a\notin M\}$.
It is easy to see that for any $R$, $\{D(a): a\in R\}$ forms a
basis of open sets on $Max(R)$. This topology is called the
Zariski topology. We use $M(I) (M(a))$ to denote the set of $M\in
Max(R)$, where $I\sub M (a\in M)$. Note that $M(I)=\bigcap_{a\in
I}M(a)$ and $M(a)=Max(R)\setminus D(a)$ (see \cite{HU}).
\\\indent For an open subset $A$ of $Max(R)$, suppose that
$M_{A}:=\{a\in R: A\sub M(a)\}$. Then we can see that $M_{A}$ is
an ideal of $R$, $M_{A}=\bigcap_{M\in A}M$ and $M(M_{A})=clA$,
where $clA$ is the cluster points of $A$ in $Max(R)$.
%----------------------------------------------
\begin{lem}\label{2.71} Let $R$ be a ring satisfying $J(R)=0$.
\begin{itemize}
\item[$(a)$] For ideals $I, J$ of $R$, $Ann(I)\sub Ann(J)$ \ifif
$intM(I)\sub intM(J).$ \item[$(b)$] For an open subset $A$ of
$Max(R)$, $M_{A}$ is an annihilator ideal. \item[$(c)$] If $I$ is
an annihilator ideal of $R$, then there exists an open subset $A$
of $Max(R)$ such that $I=M_{A}$. \item[$(d)$] For open subsets $A,
B$ of $Max(R)$, $M_{A}= M_{B}$ \ifif $clB=clA$.\item[$(e)$] $A\sub
Max(R)$ is a clopen subset \ifif there is an idempotent $e\in R$
such that $A=M(e)=D(1-e)$. \item[$(f)$] For any ideal $I$ of $R$,
$Ann(I)=M_{(\bigcup D(I))}$.
\end{itemize}
\end{lem}
\textbf{Proof.} The proof is similar to that of Lemma
\ref{2.7}.$\hfill\square$
%-----------------------------------------------------------
\begin{thm}\label{K1}
Let $R$ be a ring satisfying $J(R)=0$. Then $Max(R)$ is an
$EZ$-space \ifif for any ideal $I$ of $R$ there exists a subset
$E$ of idempotents such that $Ann(I)=Ann(E)$.
\end{thm}
\textbf{Proof.} The proof is similar to that of Theorem \ref{2.8}
$(a)\Leftrightarrow(b)$.$\hfill\square$
\\\indent Recall that a ring $R$ is called $\textit{potent}$ if idempotents
can be lifted (mod $J(R)$) and every ideal $I$ of $R$ which is
not contained in $J(R)$ contains a non-zero idempotent \cite{N}.
\begin{cor}\label{K2}
Let $R$ be a ring satisfying $J(R)=0$. If $R$ is potent, then for
any ideal $I$ of $R$ there exists a subset $E$ of idempotents such
that $Ann(I)=Ann(E)$.
\end{cor}
\textbf{Proof.} This is a consequence of \cite[Proposition 4.4]{V}
and Theorem \ref{K1}.$\hfill\square$
%-----------------------------------------
\begin{Rem}
 The converse of Corollary \ref{K2} need not be true. For
example, let $R=\Bbb Z$. Then for any nonzero ideal $I$ of $R$ we
have $Ann(I)=Ann(E)$, where $E=\{1\}$ and $J(R)=0$. But $R$ is not
a potent ring.
\end{Rem}
%-------------------------------------------------
\begin{cor}
Let $X$ be a completely regular Hausdorff space. The following
statements are equivalent.
\begin{itemize}
\item[\rm(a)] $X$ is an $EZ$-space.\item[\rm(b)] $Spec(C(X))$ is
an $EZ$-space.\item[\rm(c)] $Max(C(X))$ is an $EZ$-space.
\end{itemize}
\end{cor}
\textbf{Proof.} This is a consequence of Theorems \ref{2.5},
\ref{2.8} and \ref{K1}. $\hfill\square$

%\begin{ques}
%$(a)$ Is there an example of  space which the closure of any union
%of clopens is open but a union of clopens is not
%$C^{*}$-embedded?\\\indent $(b)$ Is there an example of space
%which disjoint unions of clopens have disjoint closure but is not
%an $EF$-space?
%\end{ques}
%{\bf ACKNOWLEDGEMENTS.}
\subsection*{Acknowledgement.}
%\begin{itemize}
The author wishes to thank Dr. M. R. Koushesh  for his
encouragement and discussion on this paper. The author also is
thankful to the referee for suggesting many improvements; in
particular, Proposition \ref{2.1} (b) and Lemma \ref{ref1} are due
to the referee.
\bibliographystyle{amsplain}

\end{document}